\documentclass[11pt,a4paper]{amsart}
\usepackage{amsfonts,amsmath,amssymb,amsthm,mathtools}
\usepackage{times}
\usepackage[noadjust]{cite}

\numberwithin{equation}{section}

\DeclareMathOperator{\Ran}{Ran}
\DeclareMathOperator{\Dom}{Dom}
\DeclareMathOperator{\spec}{spec}
\DeclareMathOperator{\dist}{dist}

\DeclarePairedDelimiter{\norm}{\lVert}{\rVert}

\newcommand{\dd}{\mathrm d}
\newcommand{\eps}{\varepsilon}

\newcommand{\N}{\mathbb{N}}
\newcommand{\EE}{\mathsf{E}}

\newcommand{\cH}{{\mathcal H}}
\newcommand{\cO}{{\mathcal O}}

\newtheorem{theorem}{Theorem}[section]{\bf}{\it}
\newtheorem{lemma}[theorem]{Lemma}{\bf}{\it}
\newtheorem{proposition}[theorem]{Proposition}{\bf}{\it}
\newtheorem{remark}[theorem]{Remark}{\it}{\rm}

\title[The subspace perturbation problem for off-diagonal perturbations]
 {Notes on the subspace perturbation problem for off-diagonal perturbations}

\author[A.\ Seelmann]{Albrecht Seelmann}
\address{A.~Seelmann, FB 08 - Institut f\"{u}r Mathematik,
Johannes Gutenberg-Universi\-t\"{a}t Mainz, Staudinger Weg 9, D-55099 Mainz, Germany}
\email{seelmann@mathematik.uni-mainz.de}

\begin{document}

\begin{abstract}
 The variation of spectral subspaces for linear self-adjoint operators under an additive bounded off-diagonal perturbation is
 studied. To this end, the optimization approach for general perturbations in [J.\ Anal.\ Math., to appear; arXiv:1310.4360 (2013)]
 is adapted. It is shown that, in contrast to the case of general perturbations, the corresponding optimization problem can not be
 reduced to a finite-dimensional problem. A suitable choice of the involved parameters provides an upper bound for the solution of
 the optimization problem. In particular, this yields a rotation bound on the subspaces that is stronger than the previously known
 one from [J.\ Reine Angew.\ Math.\ (2013), DOI:\,10.1515/crelle-2013-0099].
\end{abstract}

\maketitle

\section{Introduction}\label{sec:intro}
The present work deals with a particular case of the \emph{subspace perturbation problem} previously discussed in several recent
works such as \cite{AM13,KMM03,KMM07,MS13,Seel13}.

For the whole note let $A$ be a self-adjoint possibly unbounded operator on a separable Hilbert space $\cH$ such that the spectrum
of $A$ is separated as
\begin{equation}\label{eq:specSep}
 \spec(A) = \sigma\cup\Sigma \quad\text{ with }\quad d:=\dist(\sigma,\Sigma)>0\,.
\end{equation}
Moreover, let $V$ be a bounded self-adjoint operator on $\cH$, and suppose that $V$ is \emph{off-diagonal} with respect to the
decomposition $\cH=\Ran\EE_A(\sigma)\oplus\Ran\EE_A(\Sigma)$, that is,
\begin{equation}\label{eq:offdiag}
 \EE_A(\sigma)V\EE_A(\sigma) = 0 = \EE_A(\Sigma)V\EE_A(\Sigma)\,.
\end{equation}
Here, $\EE_A$ denotes the spectral measure for the self-adjoint operator $A$.

In this situation, it has been shown in \cite[Proposition 2.5.22]{Tre08} (see also \cite[Theorem 1.3]{KMM07} for the case of
bounded operators $A$) that
\[
 \spec(A+V) \subset \overline{\cO_{\delta_V}\bigl(\spec(A)\bigr)} \quad\text{ with }\quad
 \delta_V := \norm{V}\tan\Bigl(\frac{1}{2}\arctan2\frac{\norm{V}}{d}\Bigr)\,,
\]
where $\cO_{\delta_V}\bigl(\spec(A)\bigr)$ denotes the open $\delta_V$-neighbourhood of the spectrum of $A$. In particular, if
\begin{equation}\label{eq:gap}
 \norm{V}<\frac{\sqrt{3}}{2}\,d\,,
\end{equation}
that is, $\delta_V<d/2$, then the spectrum of $A+V$ is likewise separated into two disjoint components, each contained in the open
$d/2$-neighbourhood of $\sigma$ and $\Sigma$, respectively. This \emph{gap non-closing condition} on $\norm{V}$ is known to be
sharp, see \cite[Example 1.5]{KMM07}.

The variation of the corresponding spectral subspaces under the perturbation can be measured by the associated \emph{maximal angle}
\[
 \theta:=\arcsin\bigl(\norm{\EE_A(\sigma)-\EE_{A+V}\bigl(\cO_{d/2}(\sigma)\bigr)}\bigr)\,.
\]
Here, it is a natural question whether the condition \eqref{eq:gap} is sufficient to ensure that $\theta<\pi/2$. More specifically,
one is interested in the best possible constant $c_\text{opt-off}\in\bigl(0,\frac{\sqrt{3}}{2}\bigr]$ such that
\begin{equation}\label{eq:defcoff}
 \theta<\frac{\pi}{2} \quad\text{ whenever }\quad \norm{V}<c_\text{opt-off}\cdot d\,.
\end{equation}

It is worth mentioning that the situation looks slightly different if the two sets $\sigma$ and $\Sigma$ are additionally assumed
to be subordinated, that is, $\sup\sigma<\inf\Sigma$ or vice versa, or if one of the two sets lies in a finite gap of the other
one. In these particular cases, the gap non-closing condition on $\norm{V}$ can be relaxed considerably, and, at the same time,
these relaxed conditions are also known to ensure that the associated maximal angle $\theta$ is strictly less than $\pi/2$. A short
survey of the corresponding results can be found, e.g., in \cite{AM14}.

Under the assumption \eqref{eq:specSep} alone, however, the value of $c_\text{opt-off}$ in \eqref{eq:defcoff} is unknown. It has
been conjectured to be $\sqrt{3}/2=0{.}8660254\dots$ (see \cite{KMM07}), the same constant as in the gap non-closing condition, but
no proof for this is available yet. 

The analogous problem for general, not necessarily off-diagonal, perturbations has been discussed in \cite{KMM03,MS13,AM13,Seel13}.
In this more general setting, the (sharp) gap non-closing condition is known to be $\norm{V}<d/2$, and it is likewise an open
problem whether the best possible constant $c_\text{opt}\in\bigl(0,\frac{1}{2}\bigr]$ corresponding to \eqref{eq:defcoff} satisfies
$c_\text{opt}=1/2$. The currently best known lower bound $c_\text{opt}\ge c_\text{crit}$ with an explicit constant
$c_\text{crit}=0{.}4548\ldots$ was obtained in the author's work \cite{Seel13}. The present note ties in with the considerations
there.

As a direct consequence of the results for general perturbations, one clearly has $c_\text{opt-off}\ge c_\text{opt}$. In
particular, the maximal angle $\theta$ satisfies the bound
\begin{equation}\label{eq:sin2theta}
 \theta \le \frac{1}{2}\arcsin\Bigl(\pi\,\frac{\norm{V}}{d}\Bigr) \quad\text{ for }\quad \norm{V}\le\frac{d}{\pi}\,,
\end{equation}
proved in \cite[Corollary 2]{Seel14} (see also \cite[Corollary 4.3 and Remark 4.4]{AM13}). For off-diagonal perturbations $V$, this
bound agrees with the one obtained by combining \cite[Theorems 3.6\,(i) and 7.6]{AMM03} and \cite[Corollary 3.4]{KMM03:181}.

Making use of the off-diagonal structure of the perturbation, it was shown in \cite[Theorem 3.3]{MS13} that
\begin{equation}\label{eq:MS}
 \theta \le \frac{\pi}{2}\int_0^{\frac{\norm{V}}{d}}\frac{\dd\tau}{1-2\tau\tan\bigl(\frac{1}{2}\arctan(2\tau)\bigr)}<\frac{\pi}{2}
 \quad\text{ for }\quad \norm{V}<c_\text{off}\cdot d\,,
\end{equation}
where $c_\text{off}=0{.}6759893\ldots$ is determined by
\[
 \int_0^{c_\text{off}} \frac{\dd\tau}{1-2\tau\tan\bigl(\frac{1}{2}\arctan(2\tau)\bigr)} = 1\,.
\]
In particular, this yields the stronger lower bound $c_\text{opt-off}\ge c_\text{off}$. An earlier, slightly weaker result can be
found in \cite[Theorem 2.2]{KMM07}.

It should be noted that the result \eqref{eq:MS} was originally formulated only for the case where the operator $A$ is assumed to
be bounded, but it can easily be extended to the unbounded case. For the sake of completeness, a corresponding proof is reproduced
in Remark \ref{rem:MS} below.

In their recent survey article \cite{AM14}, Albeverio and Motovilov have stated that $c_\text{opt-off}>0{.}692834$, based on the
iteration approach from \cite{AM13} and \cite{Seel13} adapted to the case of off-diagonal perturbations.

In the present note this approach is refined. The principal result is that
\[
 c_\text{opt-off} > 0{.}6940725\,,
\]
which, together with a corresponding bound on the maximal angle, is obtained by a suitable choice of the involved parameters, see
Theorem \ref{thm:main} below. We also show that, in contrast to the case of general perturbations in \cite{Seel13}, the
optimization problem for these parameters can not be reduced to a finite-dimensional problem, see Proposition \ref{prop:noMin}
below. In fact, this optimization problem is not solved explicitly yet. Nevertheless, the result presented here is the strongest
one obtained for this problem so far.

\section{The optimization problem for off-diagonal perturbations}

Suppose that the off-diagonal perturbation $V$ is non-trivial, that is, $V\neq 0$. For $0\le t<\sqrt{3}/2$, introduce
$B_t:=A+td\cdot V/\norm{V}$ on $\Dom(B_t):=\Dom(A)$, and denote by $P_t:=\EE_{B_t}\bigl(\cO_{d/2}(\sigma)\bigr)$ the spectral
projection for $B_t$ associated with the open $d/2$-neighbourhood $\cO_{d/2}(\sigma)$ of $\sigma$.

Clearly, one has $\norm{B_t-A}=td<\sqrt{3}d/2$. As described in Section \ref{sec:intro}, this implies that the spectrum of the
perturbed operator $B_t$ is separated as
\[
 \spec(B_t) = \omega_t \cup \Omega_t
\]
with
\[
 \omega_t = \spec(B_t)\cap\overline{\cO_{\delta_t\cdot d}(\sigma)} \quad\text{ and }\quad
 \Omega_t = \spec(B_t)\cap\overline{\cO_{\delta_t\cdot d}(\Sigma)}\,,
\]
where
\[
 \delta_t := t\tan\Bigl(\frac{1}{2}\arctan 2t\Bigr)=\frac{1}{2}\sqrt{1+4t^2}-\frac{1}{2}<\frac{1}{2}\,.
\]
In particular, for $0\le t<\sqrt{3}/2$ one has $P_t=\EE_{B_t}(\omega_t)$ and
\begin{equation}\label{eq:specDist}
 \dist(\omega_t,\Omega_t) \ge (1-2\delta_t)d = \bigl(2-\sqrt{1+4t^2}\,\bigr)d\,.
\end{equation}

Let $t\in\bigl(0,\frac{\sqrt{3}}{2}\bigr)$, and let $0=t_0<\dots<t_{n+1}=t$ with $n\in\N_0$ be a finite partition of the interval
$[0,t]$. As in \cite{AM13} and \cite{Seel13}, the triangle inequality for the maximal angle (see, e.g.,
\cite[Corollary 4]{Brown93}) yields
\begin{equation}\label{eq:projTriangle}
 \arcsin\bigl(\norm{P_0-P_t}\bigr) \le \sum_{j=0}^n \arcsin\bigl(\norm{P_{t_j}-P_{t_{j+1}}}\bigr)\,.
\end{equation}
Moreover, considering $B_{t_{j+1}}=B_{t_j}+(t_{j+1}-t_j)d\cdot V/\norm{V}$ as a perturbation of $B_{t_j}$ and taking into account
the a priori bound \eqref{eq:specDist}, we observe that
\begin{equation}\label{eq:localPert}
 \frac{\norm{B_{t_{j+1}}-B_{t_j}}}{\dist(\omega_{t_j},\Omega_{t_j})} \le \frac{t_{j+1}-t_j}{1-2\delta_{t_j}} =:\lambda_j
 \quad\text{ for }\quad j=0,\dots,n\,.
\end{equation}
In particular, the bound \eqref{eq:sin2theta} for general perturbations implies that
\begin{equation}\label{eq:localSin2Theta}
 \arcsin\bigl(\norm{P_{t_j}-P_{t_{j+1}}}\bigr) \le \frac{1}{2}\arcsin(\pi\lambda_j) \quad\text{ whenever }\quad
 \lambda_j\le\frac{1}{\pi}\,.
\end{equation}

For partitions of the interval $[0,t]$ with arbitrarily small mesh size, this allows one to reproduce the bound \eqref{eq:MS} from
\cite{MS13}:

\begin{remark}[cf.\ {\cite[Section 2]{Seel13}}]\label{rem:MS}
 If the mesh size of the partition of the interval $[0,t]$ is sufficiently small, then the Riemann sum
 \[
  \sum_{j=0}^n\lambda_j = \sum_{j=0}^n \frac{t_{j+1}-t_j}{1-2\delta_{t_j}}
 \]
 is close to the integral $\int_0^t \frac{\dd\tau}{1-2\delta_\tau}$. Since, at the same time, each $\lambda_j$ is small and
 $\arcsin(x)/x\to 1$ as $x\to 0$, we conclude from \eqref{eq:projTriangle} and \eqref{eq:localSin2Theta} that
 \[
  \arcsin\bigl(\norm{P_0-P_t}\bigr) \le \frac{\pi}{2}\int_0^t \frac{\dd\tau}{1-2\delta_\tau}\,.
 \]
 Taking into account that $A+V=B_t$ for $t=\norm{V}/d$, this agrees with \eqref{eq:MS}.

 Clearly, the same reasoning for the interval $[s,t]$ can be used to show that
 \begin{equation}\label{eq:MSCont}
  \arcsin\bigl(\norm{P_s-P_t}\bigr) \le \frac{\pi}{2}\int_s^t \frac{\dd\tau}{1-2\delta_\tau} \quad\text{ for }\quad
  0\le s<t<\frac{\sqrt{3}}{2}\,.
 \end{equation}
\end{remark}

It turns out that the bound \eqref{eq:MS} on the maximal angle is stronger than \eqref{eq:sin2theta}, see Lemma
\ref{lem:offIterComp}\,(a) below. However, part (b) of the same lemma indicates that the situation changes when the estimate on the
maximal angle is iterated.

\begin{lemma}\label{lem:offIterComp}
 \hspace*{2cm}
 \begin{enumerate}
  \renewcommand{\theenumi}{\alph{enumi}}
  \item One has
        \[
         \frac{\pi}{2}\int_0^s \frac{\dd\tau}{1-2\delta_\tau} < \frac{1}{2}\arcsin(\pi s)
         \quad\text{ for }\quad 0<s\le\frac{1}{\pi}\,.
        \]
  \item For every $0<r<\sqrt{3}/2$ there is $\eps>0$ with $\eps\le(1-2\delta_r)/\pi$ and $r+\eps<\sqrt{3}/2$ such that
        \[
         \frac{1}{2}\arcsin\Bigl(\pi\,\frac{s-r}{1-2\delta_r}\Bigr) < \frac{\pi}{2}\int_r^s \frac{\dd\tau}{1-2\delta_\tau}
         \quad\text{ for }\quad r<s\le r+\eps\,.
        \]
        Moreover, the number $\eps$ can be chosen independently of $r$ from a compact subinterval of
        $\bigl(0,\frac{\sqrt{3}}{2}\bigr)$.
 \end{enumerate}

  \begin{proof}
   Let $r$ with $0\le r<\sqrt{3}/2$ be arbitrary, and define
   \[
    h_r(s) := \frac{\pi}{2}\int_r^s \frac{\dd\tau}{1-2\delta_\tau} - \frac{1}{2}\arcsin\Bigl(\pi\,\frac{s-r}{1-2\delta_r}\Bigr)\,.
   \]
   Taking into account that $1-2\delta_\tau=2-\sqrt{1+4\tau^2}$, one computes
   \[
    h_r'(s) =
    \frac{\pi}{2}\biggl( \frac{1}{2-\sqrt{1+4s^2}}-\frac{1}{\sqrt{\bigl(2-\sqrt{1+4r^2}\,\bigr){}^2-\pi^2(s-r){}^2}} \biggr)\,.
   \]
   
   For $r=0$, the inequality $h_0'(s)<0$ is equivalent to $\sqrt{1-\pi^2s^2}<2-\sqrt{1+4s^2}$, and it
   is easy to verify that the latter is valid for $0< s\le 1/\pi$. Since $h_0(0)=0$, this implies that $h_0(s)<0$ for
   $0<s\le 1/\pi$, which proves (a).
   
   Now, let $r>0$. In this case, the inequality $h_r'(s)>0$ is equivalent to
   \[
    \bigl( 2-\sqrt{1+4r^2}\, \bigr)^2 - \pi^2(s-r)^2 > \bigl( 2-\sqrt{1+4s^2}\, \bigr)^2\,,
   \]
   which, in turn, can be rewritten as
   \[
    4(s^2-r^2) \Bigl( \frac{4}{\sqrt{1+4r^2}+\sqrt{1+4s^2}} - 1 \Bigr) > \pi^2(s-r)^2\,.
   \]
   Dividing the latter inequality for $s>r$ by $s-r$ and then letting $s$ approach $r$, one arrives at the inequality
   \begin{equation}\label{eq:itEst}
    8r\cdot\Bigl( \frac{2}{\sqrt{1+4r^2}}-1 \Bigr) > 0\,,
   \end{equation}
   which is obviously valid for $0<r<\sqrt{3}/2$. Hence, by continuity, one concludes that $h_r'(s)>0$ if $s>r$ is sufficiently
   close to $r$. Since $h_r(r)=0$, this proves the first claim in (b). The second claim that $\eps$ can be chosen independently of
   $r$ from a compact subinterval of $\bigl(0,\frac{\sqrt{3}}{2}\bigr)$ follows by the same reasoning and the fact that the
   left-hand side of \eqref{eq:itEst} is bounded away from $0$ for $r$ from a compact subinterval of
   $\bigl(0,\frac{\sqrt{3}}{2}\bigr)$.
  \end{proof}%
\end{lemma}

The preceding lemma demonstrates one of the main differences between the current case of off-diagonal perturbations and the one of
general perturbations. Namely, the function $\tau\mapsto1-2\delta_\tau$ from the lower bound \eqref{eq:specDist} is not affine,
which corresponds to the fact that for $B_s=B_r+(B_s-B_r)$ with $0<r<s$ the perturbation $B_s-B_r$ does not need to be (and usually
is not) off-diagonal with respect to the decomposition $\cH=\Ran P_r\oplus\Ran P_r^\perp$. By contrast, the corresponding function
$\tau\mapsto1-2\tau$ for general perturbations is affine. In particular, for $0\le r<s<1/2$ with
$\lambda:=\frac{s-r}{1-2r}\le\frac{1}{\pi}$ one has
\[
 \frac{\pi}{2}\int_r^s \frac{\dd\tau}{1-2\tau} = \frac{\pi}{2} \int_0^\lambda \frac{\dd\tau}{1-2\tau} >
 \frac{1}{2}\arcsin(\pi\lambda)\,,
\]
cf.\ \cite[Remark 5.5]{AM13}, regardless of whether $r>0$ or $r=0$. The effect expressed by Lemma \ref{lem:offIterComp} is
therefore not present in the case of general perturbations.

In view of \eqref{eq:projTriangle}--\eqref{eq:MSCont}, Lemma \ref{lem:offIterComp} suggests to estimate $\arcsin(\norm{P_0-P_t})$
as
\begin{equation}\label{eq:optoff:Approach}
 \arcsin\bigl(\norm{P_0-P_t}\bigr) \le \frac{\pi}{2}\int_0^{\lambda_0} \frac{\dd\tau}{1-2\delta_\tau} +
 \frac{1}{2}\sum_{j=1}^n\arcsin(\pi\lambda_j)\,,
\end{equation}
provided that $\lambda_j\le1/\pi$ for $j=1,\dots,n$ and that $\lambda_0=t_1\le c_{\mathrm{off}}$ with $c_{\mathrm{off}}$ as in
\eqref{eq:MS}. The optimization problem then consists in minimizing the right-hand side of \eqref{eq:optoff:Approach} for fixed
$t\in\bigl(0,\frac{\sqrt{3}}{2}\bigr)$ over all corresponding choices of partitions of the interval $[0,t]$. This is the natural
adaption of the approach in \cite{AM13} and \cite{Seel13} to the current situation of off-diagonal perturbations $V$.

The following result shows that this optimization problem, unlike the one in \cite{Seel13}, can not be reduced to a
finite-dimensional problem. It is a direct application of Lemma \ref{lem:offIterComp}\,(b).

\begin{proposition}\label{prop:noMin}
 For fixed $t\in\bigl(0,\frac{\sqrt{3}}{2}\bigr)$, there is no finite partition of the interval $[0,t]$ which minimizes the
 right-hand side of \eqref{eq:optoff:Approach}.
 \begin{proof}
  Let $0=t_0<\dots<t_{n+1}=t$ with $n\in\N_0$ be an arbitrary partition of the interval $[0,t]$. For every $r\in(0,t_1)$ one has
  \[
   \int_0^{t_1} \frac{\dd\tau}{1-2\delta_\tau} = \int_0^r \frac{\dd\tau}{1-2\delta_\tau} +
   \int_r^{t_1} \frac{\dd\tau}{1-2\delta_\tau}\,.
  \]
  Since the number $\eps$ in Lemma \ref{lem:offIterComp}\,(b) can be chosen independently of $r$ from a compact subinterval
  of $\bigl(0,\frac{\sqrt{3}}{2}\bigr)$, we may choose $r\in(0,t_1)$ in such a way that
  $\frac{t_1-r}{1-2\delta_r}\le\frac{1}{\pi}$ and
  \[
   \frac{1}{2}\arcsin\Bigl(\pi\,\frac{t_1-r}{1-2\delta_r}\Bigr) < \frac{\pi}{2}\int_r^{t_1} \frac{\dd\tau}{1-2\delta_\tau}\,.
  \]
  The refined partition $0=t_0<r<t_1<\dots<t_{n+1}=t$ then leads to a right-hand side in \eqref{eq:optoff:Approach} which is
  strictly less than the one corresponding to the original partition of the interval $[0,t]$.
 \end{proof}%
\end{proposition}

\begin{remark}
 Iterating the argument in the proof of Proposition \ref{prop:noMin} shows that, from the point of view of minimizing the
 right-hand side of \eqref{eq:optoff:Approach}, one may always assume that $\lambda_0=0$. In other words, the first summand of the
 right-hand side of \eqref{eq:optoff:Approach} can be replaced by $\frac{1}{2}\arcsin(\pi\lambda_0)$, provided that
 $\lambda_0\le1/\pi$, without affecting the optimization result. In fact, the latter has been considered in \cite{AM14}. However,
 when considering specific finite partitions of the interval $[0,t]$, Lemma \ref{lem:offIterComp}\,(a) shows that the current
 approach \eqref{eq:optoff:Approach} is more suitable.
\end{remark}

Another difficulty in the problem of minimizing the right-hand side of \eqref{eq:optoff:Approach} arises by the fact that, given a
partition of the interval $[0,t]$, an efficient explicit representation of $t$ in terms of the corresponding parameters $\lambda_j$
is not at hand. In contrast to the case of general perturbations (cf.\ \cite[Section 3]{Seel13}), it is thus unclear how to
determine the critical points for the reduced finite-dimensional optimization problems associated with a fixed number of supporting
points in the partitions.

In fact, the problem of minimizing the right-hand side of \eqref{eq:optoff:Approach} is not solved explicitly yet. So far, the
author can only guess a choice of the parameters $\lambda_j$ guaranteeing that $c_\text{opt-off} > 0{.}694$. In view of Proposition
\ref{prop:noMin}, this guess seems to be a reasonable compromise between the complexity of the choice of the parameters and the
strength of the result:

Let $n=4$. Choose $\lambda_0\in(0,c_\text{off})$ such that
\[
 \frac{\pi}{2}\int_0^{\lambda_0} \frac{\dd\tau}{1-2\delta_\tau} = \frac{1}{3}
\]
and $\lambda_j\in\bigl(0,\frac{1}{\pi}\bigr]$, $j=1,\dots,4$, such that $2\arcsin(\pi\lambda_j)=\frac{\pi}{2}-\frac{1}{3}$, that
is,
\[
 \lambda_j = \frac{1}{\pi}\sin\Bigl(\frac{3\pi-2}{12}\Bigr) = 0{.}1846204\dots \quad\text{ for }\quad j=1,\dots,4\,.
\]
For this choice of $n$ and $\lambda_j$ the right-hand side of \eqref{eq:optoff:Approach} equals $\pi/2$.

A numerical calculation gives
\[
 \tau_1 := \lambda_0 > 0{.}2062031\,.
\]
Upon observing that the mapping $\bigl[0,\frac{\sqrt{3}}{2}\bigr]\ni \tau\mapsto \tau+\lambda_1(1-2\delta_\tau)$ is strictly
increasing, it is then easy to verify that
\[
 \tau_2 := \tau_1 + \lambda_1(1-2\delta_{\tau_1}) > 0{.}3757396\,.
\]
In the same way, one has
\[
 \tau_3 := \tau_2 + \lambda_2(1-2\delta_{\tau_2})> 0{.}5140409 \,,\
 \tau_4:=\tau_3 + \lambda_3(1-2\delta_{\tau_3}) > 0{.}6184976\,,
\]
and
\begin{equation}\label{eq:ccrit}
 c_\text{off}^*:=\tau_5 := \tau_4 + \lambda_4(1-2\delta_{\tau_4}) > 0{.}6940725\,.
\end{equation}
Finally, consider the piecewise defined function $N_\text{off}^*\colon\bigl[0,c_\text{off}^*]\to\bigl[0,\frac{\pi}{2}\bigr]$
with
\begin{equation}\label{eq:defNoff}
 N_\text{off}^*(t) := \begin{cases}
                       \frac{\pi}{2}\int_0^t \frac{\dd\tau}{1-2\delta_\tau}\,, & 0\le t \le \tau_1\\[0.1cm]
                       \frac{1}{3} + (j-1)\,\frac{3\pi-2}{24} +
                        \frac{1}{2}\arcsin\bigl(\pi\,\frac{t-\tau_j}{1-2\delta_{\tau_j}}\bigr)\,,
                        & \tau_j < t \le \tau_{j+1}\,.
                      \end{cases}
\end{equation}
Clearly, the function $N_\text{off}^*$ is strictly increasing with $N_\text{off}^*(c_\text{off}^*)=\pi/2$, continuous on
$[0,c_\text{off}^*]$, and continuously differentiable on $(0,c_\text{off}^*)\setminus\{\tau_2,\tau_3,\tau_4\}$.

We now use $\tau_1,\dots,\tau_4$ as supporting points for the partitions of the interval $[0,t]$. More precisely, using the
partition $0<\dots<\tau_j<t$ for $\tau_j< t\le\tau_{j+1}$ and the trivial one $0<t$ for $t\le\tau_1$, it follows from
\eqref{eq:optoff:Approach} and \eqref{eq:defNoff} that
\[
 \arcsin\bigl(\norm{P_0-P_t}\bigr) \le N_\text{off}^*(t)<\frac{\pi}{2} \quad\text{ for }\quad 0\le t<c_\text{off}^*\,. 
\]

Taking into account that $B_t=A+V$ with $t=\norm{V}/d$, the preceding considerations now summarize to the following theorem, the
main result in this note.

\begin{theorem}\label{thm:main}
 Let $A$ be a self-adjoint operator on a separable Hilbert space $\cH$ with spectrum separated as in \eqref{eq:specSep}, and let
 $V$ be a bounded self-adjoint operator on $\cH$ which is off-diagonal with respect to the decomposition
 $\cH=\EE_A(\sigma)\oplus\EE_A(\Sigma)$, that is,
 \[
  \EE_A(\sigma)V\EE_A(\sigma) = 0 = \EE_A(\Sigma)V\EE_A(\Sigma)\,.
 \]
 If $V$ satisfies
 \[
  \norm{V}<c_\mathrm{off}^*\cdot d
 \]
 with $c_\mathrm{off}^*$ as in \eqref{eq:ccrit}, then
 \[
  \arcsin\bigl(\norm{\EE_A(\sigma)-\EE_{A+V}\bigl(\cO_{d/2}(\sigma)\bigr)}\bigr) \le N_\mathrm{off}^*\Bigl(\frac{\norm{V}}{d}\Bigr)
  < \frac{\pi}{2}\,,
 \]
 where the function $N_\mathrm{off}^*$ is given by \eqref{eq:defNoff}.
\end{theorem}

It is a direct consequence of Theorem \ref{thm:main} that the best possible constant $c_\text{opt-off}$ in \eqref{eq:defcoff}
satisfies the lower bound
\begin{equation}\label{eq:copt}
 c_\text{opt-off} > c_\text{off}^* > 0{.}6940725\,,
\end{equation}
where the fact that the first inequality in \eqref{eq:copt} is strict is due to Proposition \ref{prop:noMin}. Furthermore,
numerical evaluations suggest that the corresponding bound on the maximal angle between the subspaces $\Ran\EE_A(\sigma)$ and
$\Ran\EE_{A+V}\bigl(\cO_{d/2}(\sigma)\bigr)$ is indeed stronger than the one given by \eqref{eq:MS}, that is,
\[
 N_\text{off}^*(t) < \frac{\pi}{2} \int_0^t \frac{\dd\tau}{1-2\delta_\tau} \quad\text{ for }\quad \tau_1<t\le c_\text{off}\,.
\]

\section*{Acknowledgements}
The material presented in this work is part of the author's Ph.D.\ thesis \cite{SeelDiss}. The author is grateful to his
Ph.D.\ advisor Vadim Kostrykin for fruitful discussions. He would also like to thank Christoph Uebersohn for helpful remarks on the
manuscript.


\end{document}